\documentclass[proceedings,submission%
]{dmtcs}

\usepackage[latin1]{inputenc}
\usepackage{subfigure}
\usepackage{amsmath, amssymb}

\newtheorem{thm}{Theorem}

\newtheorem{cor}{Corollary}

\newcommand{\U}{\mathcal{U}}
\newcommand{\Cat}{\operatorname{Cat}}
\newcommand{\C}{\mathbb{C}}
\newcommand{\Z}{\mathcal{Z}}

\newcommand{\RC}{\operatorname{RC}}

\renewcommand{\a}{\tilde{a}}
\newcommand{\M}{\mathcal{M}}

\author{Sho Matsumoto\addressmark{1}\and Jonathan Novak\addressmark{2}}
\title[Primitive Factorizations, Jucys-Murphy Elements, and Matrix Models]{Primitive Factorizations,
Jucys-Murphy Elements, and Matrix Models}
\address{\addressmark{1}Graduate School of Mathematics, Nagoya University, Furocho, Chikusa-ku, Nagoya, 464-8602, Japan\\
\addressmark{2}Department of Combinatorics \& Optimization, University of Waterloo, Waterloo,
Canada}

\begin{document}
\maketitle
\begin{abstract}
\paragraph{Abstract.}
	A factorization of a permutation into transpositions is called ``primitive'' if its factors
	are weakly ordered.
	We discuss the problem of enumerating primitive factorizations of permutations,
	and its place in the hierarchy of previously studied factorization problems.  Several
	formulas enumerating minimal primitive and possibly non-minimal primitive
	factorizations are presented, and interesting connections with Jucys-Murphy elements,
	symmetric group characters, and matrix models are described.

\paragraph{R\'esum\'e.}
	Une factorisation en transpositions d'une  permutation est dite
	``primitive'' si ses facteurs sont ordonn\'es.  Nous discutons du
	probl\`eme de l'\'enum\'eration des factorisations primitives de
	permutations, et de sa place dans la hi\'erarchie des probl\`emes de
	factorisation pr\'ec\'edemment \'etudi\'es.  Nous pr\'esentons plusieurs
	formules \'enum\'erant certaines classes de factorisations primitives,
	et nous soulignons des connexions int\'eressantes avec les \'el\'ements
	Jucys-Murphy, les caract\'eres des groupes sym\'etriques, et les
	mod\`eles de matrices.


\end{abstract}

\section{Introduction}
The problem of counting the number of ways in which a given permutation can be 
factored into a given number of transpositions is of perennial interest in algebraic
combinatorics.  Usually, one considers this problem in the presence of constraints
on the factors, e.g. that they should be transpositions of a certain type, 
or should collectively generate a certain group, etc.
By varying these constraints, one obtains enumeration problems which enjoy surprising 
connections with other branches of mathematics.

The earliest enumerative study of transposition factorizations was carried out by 
Hurwitz \cite{Hurwitz} in the nineteenth century.  Motivated by a problem from enumerative
algebraic geometry, namely the counting of almost simple ramified covers of the sphere by 
other Riemann surfaces, Hurwitz published an explicit 
formula for the number of \emph{minimal transitive factorizations} of an arbitrary permutation into
transpositions.  There are two constraints in the Hurwitz factorization problem: ``minimality,'' 
which requires that the number of transpositions used should be 
as small as possible, and ``transitivity,'' which requires that the factors should act transitively on 
the points $\{1,\dots,n\}.$  Hurwitz's formula for the number of minimal transitive factorizations
of a permutation $\pi \in S(n)$ of cycle type $\mu=(\mu_1,\dots,\mu_\ell)\vdash n$ is
	\begin{equation}
		\label{Hurwitz}
		(n+\ell-2)!n^{\ell-3}\prod_{i=1}^{\ell} \frac{\mu_i^{\mu_i}}{(\mu_i-1)!}.
	\end{equation}
As a particularly beautiful special case, Hurwitz's formula yields
that the number of factorizations of a full cycle in $S(n)$ into $n-1$ transpositions is
$n^{n-2},$ the number of trees on $n$ labelled vertices.  This case of Hurwitz's formula 
was independently rediscovered and popularized by D\'enes \cite{Denes}.  
The general Hurwitz formula
was independently rediscovered by Goulden and Jackson \cite{GJ:transitive},
to whom the first rigorous proof is due.

A key feature of the Hurwitz factorization problem is centrality: the number of minimal transitive
factorizations of $\pi$ depends only on the cycle type of $\pi.$  This remains true for transitive
factorizations of arbitrary length.  A different choice of constraints leading to a non-central factorization
problem was considered by Stanley \cite{Stanley:reduced}, who initiated the study of what he 
termed \emph{reduced decompositions}.  These are minimal factorizations in which the
transpositions allowed to be used as factors are the Coxeter generators $(s,s+1).$
Reduced decompositions are also referred to as \emph{sorting networks} because of
their relation to the bubblesort algorithm familiar to computer scientists.  They could also be called 
\emph{minimal Coxeter factorizations}.
The enumeration of reduced decompositions is complicated by its non-centrality, and has
spawned its own extensive literature, see \cite{Garsia} for a beautiful introduction.  
The asymptotic behaviour of random reduced decompositions
is the subject of an intriguing set of conjectures due to Angel et al. \cite{AHRV}.

The Coxeter factorization problem naturally fits into a wider class
of constrained factorization problems, in which the factors are chosen from a specified set of transpositions which generate $S(n).$  A second example from
from this class was considered by Pak \cite{Pak}, who initiated the study of \emph{star factorizations}.  These are factorizations in which the transpositions allowed to be used as factors have the 
form $(1*).$  For example, the unique minimal star factorization of $(123) \in S(3)$ is 
$(123)=(13)(12).$  Thanks to recent work of Irving and Rattan \cite{IR}, Goulden and Jackson \cite{GJ:star},
and F\'eray \cite{Feray}, the combinatorics of star factorizations is now completely understood.

Recently, Gewurz and Merola \cite{GM} posed the problem of enumerating 
transposition factorizations under a constraint on the order of the factors.  A factorization
	\begin{equation}
		\label{factorization}
		\pi = (s_1,t_1) \dots (s_k,t_k) \qquad (s_i<t_i)
	\end{equation}
of $\pi \in S(n)$ into a product of $k$ transpositions is called \emph{primitive} if
	\begin{equation}
		2 \leq t_1 \leq \dots \leq t_k \leq n,
	\end{equation}
i.e. if its factors appear in weakly increasing order with respect to the larger element in each.  For
example, $(123) \in S(3)$ can be factored into a product of two transpositions
in three ways,
	\begin{equation}
		(123)=(12)(23)=(23)(13)=(13)(12),
	\end{equation}
but only the first two of these factorizations are primitive.
Gewurz and Merola obtain the interesting result that the number of primitive factorizations
of the cycle $(12\dots n) \in S(n)$ into $n-1$ transpositions (i.e. the number
of minimal primitive factorizations) is the Catalan number
	\begin{equation}
		\Cat_{n-1} = \frac{1}{n}{2n-2 \choose n-1}.
	\end{equation}
This result should be considered in tandem with Hurwitz's $n^{n-2}$-count of 
unrestricted factorizations of $(12\dots n).$

In this extended abstract prepared for FPSAC 2010, we will give 
an overview of the authors' ongoing work
on the enumeration of primitive factorizations.  First we study minimal primitive 
factorizations, and obtain an analogue of Hurwitz's formula \eqref{Hurwitz}, i.e. an explicit formula
which counts minimal primitive factorizations of an arbitrary permutation (Theorem 
\ref{thm:minimalType} and Corollary \ref{cor:minimalTotal} below).  Then we present
a link between the primitive factorization problem and Jucys-Murphy elements.  This connection
explains the centrality of the primitive factorization problem, and allows us to use 
character theory to enumerate primitive factorizations of a full cycle into any number
of transpositions (Theorem \ref{thm:centralFactorial} below).  Finally, we discuss a surprising
connection between the primitive factorization problem and the theory of matrix models:
generating functions enumerating primitive factorizations may be expressed as integrals
over groups of unitary matrices against the Haar measure (Theorem \ref{thm:matrixIntegral} below).
It turns out that these integrals are of independent interest and have a long history in 
mathematical physics.

\section{Minimal Primitive Factorizations}
Any primitive factorization of $\pi$ into $k$ transpositions has the form
	\begin{equation}
		\pi = \underbrace{(*2) \dots (*2)}_{a_2}  \underbrace{(*3) \dots (*3)}_{a_3} \dots
			 \underbrace{(*n) \dots (*n)}_{a_n},
	\end{equation}
where $(a_2,\dots,a_n)$ is a weak $(n-1)$-part composition of $k.$
We will say that the above factorization is of type $\lambda \vdash k$ if the 
frequencies $a_2,a_3,\dots,a_n$ coincide with the parts of $\lambda$ after reordering.
For example, there are three primitive factorizations of $(1234) \in S(4)$ of type $(2,1),$ namely 
	\begin{equation}
		(1234)=\underbrace{(23)(13)}_2 \underbrace{(34)}_1=
		\underbrace{(12)}_1\underbrace{(34)(24)}_2
		=\underbrace{(23)}_1 \underbrace{(34)(14)}_2.
	\end{equation}
	
Let us now enumerate minimal primitive factorizations by type.  Let $\mathfrak{E}(k)$ denote
the set of all weakly increasing sequences $i_1 \leq \dots \leq i_k$ of $k$ positive integers 
such that $i_p \geq p$ for $1 \leq p \leq k-1$ and $i_k=k.$  It is not difficult to show that
	\begin{equation}	
		|\mathfrak{E}(k)|=\Cat_k.
	\end{equation}
Given a partition $\lambda \vdash k,$ one may introduce a refinement $\RC(\lambda)$
of the Catalan number by declaring $\RC(\lambda)$ to be the number of sequences in 
$\mathfrak{E}(k)$ of type $\lambda.$  Then, by definition,
	\begin{equation}
		\sum_{\lambda \vdash k} \RC(\lambda)=\Cat_k.
	\end{equation}
These refined Catalan numbers have previously been studied by Haiman \cite{Haiman}
and Stanley \cite{Stanley:park} in connection with parking functions, and are known explicitly:
	\begin{equation}
		\RC(\lambda)= \frac{|\lambda|!}{(|\lambda|-\ell(\lambda)+1)!\prod_{i \geq 1} m_i(\lambda)!},
	\end{equation}
where $m_i(\lambda)$ is the multiplicity of $i$ in $\lambda.$  Finally, given a pair of 
partitions $\lambda,\mu \vdash k$, introduce the set of sequences of partitions
	\begin{equation}
		\mathfrak{R}(\lambda,\mu) = \{(\lambda^{(1)},\dots,\lambda^{\ell(\mu)}) | 
		\lambda^{(i)} \vdash \mu_i, \quad \lambda^{(1)} \cup \dots \cup \lambda^{(\ell(\mu))}=
		\lambda\}.
	\end{equation}
Thus sequences in $\mathfrak{R}(\lambda,\mu)$ are obtained by breaking parts of $\mu$
in such a way that, after sorting, one obtains $\lambda.$
If $\mathfrak{R}(\lambda,\mu)$ is non-empty, then $\lambda$ is said to be a refinement
of $\mu.$

\begin{thm}
	\label{thm:minimalType}
	Let $\pi \in S(n)$ be a permutation of reduced\footnote{Recall that the 
	reduced cycle type of $\pi$ is the partition obtained by subtracting one from the 
	length of each of its cycles.  Thus the size of the reduced cycle
	type of $\pi$ is the length of a minimal factorization of $\pi$ into transpositions.} 
	cycle type $\mu \vdash k$, and 
	let $\lambda$ be another partition of $k.$  The number of primitive factorizations
	of $\pi$ of type $\lambda$ is 
		$$\sum_{(\lambda^{(1)},\dots,\lambda^{(\ell(\mu))}) \in \mathfrak{R}(\lambda,\mu)}
		\RC(\lambda^{(1)}) \dots \RC(\lambda^{(\ell(\mu))}).$$
\end{thm}

As an example, consider the permutation 
	$$\pi=(123456)(78910) \in S(10).$$
This permutation has reduced cycle type $\mu=(5,3),$ so the length of a minimal 
primitive factorization of $\pi$ is eight.  By Theorem \ref{thm:minimalType}, the number
of minimal primitive factorizations of $\pi$ of type $\lambda=(3,2,2,1)$ is
	$$\RC(3,2)\RC(2,1)+\RC(2,2,1)\RC(3)=5 \cdot 3 + 10 \cdot 1 = 25.$$
	
The proof of Theorem \ref{thm:minimalType} is bijective, and we refer the reader to 
our full-length article \cite{MN} for details.  As a corollary of this result, we obtain 
an elegant formula which counts the total number of minimal primitive factorizations of an arbitrary 
permutation.  Gewurz and Merola's result is recovered as the case $\mu=(n)$ of this corollary.

\begin{cor}
	\label{cor:minimalTotal}
	Let $\pi \in S(n)$ be a permutation of non-reduced cycle type $\mu \vdash n.$  
	The total number of primitive factorizations of $\pi$ into $n-\ell(\mu)$ transpositions
	is 
		$$\prod_{i=1}^{\ell(\mu)} \Cat_{\mu_i-1}.$$
\end{cor}

\section{Primitive Factorizations and Jucys-Murphy Elements}
We will now give an algebraic explanation of the fact that the primitive factorization
problem is central.  Let $\C[S(n)]$ denote the group algebra of the symmetric group,
and $\mathcal{Z}(n)$ its center.  The center is a commutative algebra with canonical
basis $\{C_{\mu} : \mu \vdash n\}$ consisting of the conjugacy classes of $S(n);$ for
this reason we call $\mathcal{Z}(n)$ the \emph{class algebra}.

Let $\Lambda$ denote the algebra of symmetric functions over $\C.$  We define a specialization
$\Lambda \rightarrow \Z(n)$ as follows.  For $k \geq 1,$ put
	\begin{equation}
		J_k := \sum \text{transpositions in }S(k) - \sum \text{transpositions in }S(k-1)
		=(1,k) + \dots + (k-1,k).
	\end{equation}
Thus $J_1,\dots,J_n \in \C[S(n)],$ with $J_1=0.$  The elements so defined are called 
\emph{Jucys-Murphy elements}.  They were introduced independently by Jucys \cite{Jucys} and Murphy \cite{Murphy}.  These simple elements have many remarkable properties, some of which
we will make use of here.  Diverse applications of the Jucys-Murphy elements are found in 
the work of Okounkov \cite{Okounkov96,Okounkov00} and Okounkov and Vershik \cite{OV}.

Although $\{J_1,\dots,J_n\} \not \subseteq \Z(n)$ for $n \geq 3,$ the JM elements
do belong to the Gelfand-Zetlin subalgebra of $\C[S(n)].$  This is the maximal commutative
subalgebra of $\C[S(n)]$ generated by the class algebras $\Z(1),\dots,\Z(n),$ where 
$\Z(1),\dots,\Z(n-1)$ are embedded in $\C[S(n)]$ in the canonical way.  This is clear, since
$J_k$ is by definition the difference of an element of $\Z(k)$ and an element of $\Z(k-1).$
Consequently,
the JM elements commute with one another, and we may define the alphabet 
$\Xi_n=\{\{J_1,\dots,J_n,0,0,\dots\}\}$ and evaluate symmetric functions on this 
alphabet.  It is a remarkable result of Jucys that, for any $f \in \Lambda,$
	\begin{equation}
		\label{Jucys}
		f(\Xi_n) \in \Z(n).
	\end{equation}
Thus symmetric functions of JM elements are central,
and we have the \emph{JM specialization} $\Lambda \rightarrow \Z(n).$

The JM specialization gives a very clean proof of the centrality of the primitive factorization
problem.  Let $h_k \in \Lambda$ denote the complete homogoneous symmetric function
of degree $k,$ i.e.
	\begin{equation}
		h_k=\sum_{1 \leq i_1 \leq \dots \leq i_k} x_{i_1} \dots x_{i_k},
	\end{equation}
the sum of all degree $k$ monomials in the variables $x_1,x_2,\dots.$
Then, by Jucys' result, $h_k(\Xi_n) \in \Z(n)$ and we may write
	\begin{equation}
		h_k(\Xi_n) = \sum_{\mu \vdash n} a_{k,\mu} C_{\mu}
	\end{equation}
for some coefficients $a_{k,\mu}.$  On the other hand, we have that
	\begin{equation}
		\begin{split}
		h_k(\Xi_n) &= \sum_{2 \leq t_1 \leq \dots \leq t_k \leq n} J_{t_1} \dots J_{t_k} \\
		&= \sum_{2 \leq t_1 \leq \dots \leq t_k \leq n} \sum_{s_1 < t_1} (s_1,t_1)
		\dots \sum_{s_k < t_k} (s_k,t_k) \\
		&=\sum_{\pi \in S(n)} \#\{\text{primitive factorizations of $\pi$ into $k$ transpositions}\}
		\pi.
		\end{split}
	\end{equation}
Thus, for any permutation $\pi \in S(n)$ of cycle type $\mu,$ we have that
	\begin{equation}	
		\label{primitiveCount}
		a_{k,\mu} = \#\{\text{primitive factorizations of $\pi$ into $k$ transpositions}\}.
	\end{equation}
In other words, the combinatorial problem of counting primitive factorizations is equivalent
to the algebraic problem of resolving $h_k(\Xi_n)$ with respect to the canonical basis 
of the class algebra $\Z(n).$  More generally, for any partition $\lambda$ one may consider
the resolution
	\begin{equation}
		m_{\lambda}(\Xi_n) = \sum_{\mu \vdash n} b_{\lambda\mu} C_{\mu},
	\end{equation}
where $m_{\lambda}$ is the monomial symmetric function of type $\lambda.$  We then
have, for any $\pi \in C_{\mu},$
	\begin{equation}
		b_{\lambda\mu} = \#\{\text{primitive factorizations of $\pi$ of type $\lambda$}\}.
	\end{equation}
	
\section{The Case of a Single Cycle}
The above algebraic encoding of the primitive factorization problem
allows us to enumerate primitive factorizations of a full cycle $\pi \in C_{(n)}.$
This is thanks to the remarkable properties of JM elements in irreducible representations
of $\C[S(n)],$ which amount to the fact that, while it is difficult to compute the the coordinates
of $f(\Xi_n)$ with respect to the class basis of $\Z(n),$ it is easy to compute its coordinates
with respect to the character basis of $\Z(n).$

The following remarkable result is due to Jucys \cite{Jucys}, see \cite{Okounkov96} for a complete
proof.
Given a Young diagram $\lambda$ and a cell $\Box$ in $\lambda,$ recall that the 
\emph{content} of $\lambda$ is defined to be the column index of $\lambda$ less its
row index.  Let us associate to each Young diagram $\lambda$ its alphabet 
	\begin{equation}
		A_{\lambda} = \{\{c(\Box) : \Box \in \lambda\}\}
	\end{equation}
of contents.  Let $H_{\lambda}$ denote the product of all hook-lengths of $\lambda.$
Let $\{\chi^{\lambda} : \lambda \vdash n\}$ be the characters of the irreducible
representations of $S(n),$ which form a basis of $\Z(n).$
  Then, for any symmetric function $f \in \Lambda,$ the character expansion
of $f(\Xi_n)$ is 
	\begin{equation}
		\label{characterExpansion}
		f(\Xi_n)=\sum_{\lambda \vdash n} \frac{f(A_{\lambda})}{H_{\lambda}}\chi^{\lambda}.
	\end{equation}
Thus, at the combinatorial level, the character expansion of $f(\Xi_n)$ is implemented 
by the substitution rule $\Xi_n \rightarrow A_{\lambda}.$

Let us use the above character expansion result to enumerate primitive factorizations
of a full cycle $\pi \in C_{(n)}$ of any given length.  We already know that the number
of minimal primitive factorizations of $(12 \dots n),$ i.e. those consisting of $n-1$ transpositions,
is the Catalan number $\Cat_{n-1}.$  We will now solve the problem when the number of
transpositions used is allowed to be arbitrary.  As we will see momentarily, this problem
has an unexpectedly simple and beautiful solution.

Let $z$ be an indeterminate, and form the generating function
	\begin{equation}
		\Phi(z;n) = \sum_{k \geq 0} h_k(\Xi_n) z^k.
	\end{equation}
This generating function is an element of the algebra $\Z(n)[[z]]$ of single-variable formal power series
with coefficients in the class algebra $\Z(n).$  By Jucys' character expansion result, we have
	\begin{equation}
		\begin{split}
			\Phi(z;n) &= \sum_{k \geq 0} \bigg{(} \sum_{\lambda \vdash n} 
				\frac{h_k(A_{\lambda})}{H_{\lambda}} \chi^{\lambda} \bigg{)} z^k \\
				&= \sum_{\lambda \vdash n}  \bigg{(} \sum_{k \geq 0} h_k(A_{\lambda})z^k
				 \bigg{)} \frac{\chi^{\lambda}}{H_{\lambda}} \\
				 &=  \sum_{\lambda \vdash n} \frac{\chi^{\lambda}}{H_{\lambda}
				 \prod_{\Box \in \lambda}(1- c(\Box)z)}
		\end{split}
	\end{equation}
where the last line follows from the generating function
	\begin{equation}
		\sum_{k \geq 0} h_k(x_1,\dots,x_n)z^k = \prod_{i=1}^n \frac{1}{1-x_i z}
	\end{equation}
for the elementary symmetric functions.  Note that this computation shows that the generating
function $\Phi(z;n)$ is actually a \emph{rational} function over $\Z(n).$ 

Now, given $\mu \vdash n,$ in order to obtain the generating function
	\begin{equation}
		\Phi_{\mu}(z) = \sum_{k \geq 0} a_{k,\mu} z^k
	\end{equation}
we simply take the corresponding traces of the conjugacy class $C_{\mu}$ in each irreducible
representation, to obtain a rational function over $\C:$
	\begin{equation}
		\Phi_{\mu}(z) = \sum_{\lambda \vdash n} \frac{\chi^{\lambda}(C_{\mu})}{H_{\lambda}
				 \prod_{\Box \in \lambda} (1-c(\Box)z)}.
	\end{equation}
Up until this point, the partition $\mu \vdash n$ has been generic, but now we restrict to
the special case $\mu=(n),$ the partition of $n$ with a single part.
A classical result from representation theory informs us that the trace of $C_{(n)}$ in an 
irreducible representation can only be non-zero in ``hook'' representations:
	\begin{equation}
		\chi^{\lambda}(C_{(n)})= \begin{cases}
			(-1)^{r}, \text{ if } \lambda=(n-r,1^r) \\
			0, \text{ otherwise }
			\end{cases}.
	\end{equation}
Now, the content alphabet of a hook diagram may be obtained immediately,
	\begin{equation}
		A_{(n-r,1^r)}=\{0,1,\dots,n-r-1\} \sqcup \{-1,\dots,-r\}.
	\end{equation}
so that
	\begin{equation}
		\Phi_{(n)}(z)=\sum_{r=0}^{n-1} \frac{(-1)^r}{H_{(n-r,1^r)}\prod_{i=1}^{n-r-1}(1-iz)
		\prod_{j=1}^r (1+jz)}.
	\end{equation}
For example, if $n=4,$ this is a rational function of the form
	\begin{equation}
		\begin{split}
		\Phi_{(4)}(z) = &\frac{\text{const.}}{(1-z)(1-2z)(1-3z)}
				+\frac{\text{const.}}{(1-z)(1-2z)(1+z)} \\
				+&\frac{\text{const.}}{(1-z)(1+z)(1+2z)}
				+\frac{\text{const.}}{(1+z)(1+2z)(1+3z)}.
		\end{split}
	\end{equation}
Thus, as an irreducible rational function, $\Phi_{(n)}(z)$ has the form
	\begin{equation}
		\Phi_{(n)}(z) = \frac{\sum_{i=0}^{n-1} c_i z^i}{\prod_{i=1}^{n-1}(1-i^2z^2)}
	\end{equation}
where $c_0,\dots,c_{n-1} \in \C$ are some constants to be determined momentarily.

Before finding the above coefficients, let us consider the generating function
	\begin{equation}
		\frac{1}{\prod_{i=1}^n (1-i^2u)}=\sum_{g \geq 0} h_g(1^2,\dots,n^2)u^g.
	\end{equation}
The coefficients in this generating function are complete symmetric functions evaluated
on the alphabet $\{1^2,\dots,n^2\}$ of square integers.  Reason dictates that they ought to 
be close relatives of the Stirling numbers
	\begin{equation}
		S(n+g,n)=h_g(1,\dots,n).
	\end{equation}
The Stirling number $S(a,b)$ has the following combinatorial interpretation: it counts 
the number of partitions
	\begin{equation}
		\{1,\dots,a\} = V_1 \sqcup \dots \sqcup V_b
	\end{equation}
of an $a$-element set into $b$ disjoint non-empty subsets.  Stirling numbers are given
by the explicit formula
	\begin{equation}
		S(a,b) = \sum_{j=0}^b (-1)^{b-j} \frac{j^a}{j!(b-j)!}.
	\end{equation}
The numbers 
	\begin{equation}
		T(n+g,n)=h_g(1^2,\dots,n^2)
	\end{equation}
are known as \emph{central factorial numbers}.  The central factorial numbers were 
studied classically by Carlitz and Riordan, see \cite[Exercise 5.8]{Stanley} for references.
They have the following combinatorial interpretation: $T(a,b)$ counts the number of partitions
	\begin{equation}
		\{1,1',\dots,a,a'\} = V_1 \sqcup \dots \sqcup V_b
	\end{equation}
of a set of $a$ marked and $a$ unmarked points into $b$ disjoint non-empty subsets such 
that\footnote{B\'alint Vir\'ag gave a colourful description of this condition, which is actually 
quite a useful mnemonic: ``the most important guy gets to bring his wife.''}, for each block $V_j,$ if $i$ is the least integer such that either $i$ or $i'$ appears in 
$V_j,$ then $\{i,i'\} \subseteq V_j.$  Central factorial numbers are given by the explicit formula
	\begin{equation}
		T(a,b)=2\sum_{j=0}^b (-1)^{b-j} \frac{j^{2a}}{(b-j)!(b+j)!}.
	\end{equation}

Now let us determine the unknown constants $c_0,\dots,c_{n-1}.$  By the above discussion,
the generating function $\Phi_{(n)}(z)$ has the form
	\begin{equation}
		\Phi_{(n)}(z)=(c_0+c_1z+\dots+c_{n-1}z^{n-1}) \sum_{g \geq 0}T(n-1+g,n-1)z^{2g}.
	\end{equation}
On the other hand, by the results of the previous sections,
	\begin{equation}
		\begin{split}
		\Phi_{(n)}(z) &=\sum_{k \geq 0} a_{k,(n)}z^k \\
		&= \sum_{g \geq 0} a_{n-1+2g,(n)} z^{n-1+2g}, \text{ since every permutation is either 
				even or odd}, \\
		&= \Cat_{n-1}z^{n-1} + a_{n+1,(n)}z^{n+1} + \dots, \text{ since }a_{n-1,(n)}=\Cat_{n-1}.
		\end{split}
	\end{equation}
Consequently, we must have $c_0= \dots = c_{n-2}=0, c_{n-1}=\Cat_{n-1},$ and we have
proved the following result.

\begin{thm}
	\label{thm:centralFactorial}
	For any $g \geq 0,$ the number of primitive factorizations of $(12\dots n) \in S(n)$ into $n-1+2g$
	transpositions is 
		$$\Cat_{n-1} \cdot T(n-1+g,n-1),$$
	where $T(a,b)$ denotes the Carlitz-Riordan central factorial number.  Equivalently,
	we have the generating function
		$$\Phi_{(n)}(z) = \frac{\Cat_{n-1}z^{n-1}}{(1-1^2z^2) \dots (1-(n-1)^2z^2)}.$$
\end{thm}

\section{Primitive Factorizations and Matrix Models}
Finally, we come to what is perhaps the most intriguing aspect of the primitive factorization
problem: its connection with matrix models.  The theory of matrix models has its origins in 
an area of mathematical physics known as quantum field theory, see \cite{Etingof} for 
a solid introduction.  
For our purposes, the following grossly oversimplified description of matrix model theory suffices:
	\begin{enumerate}
	
		\item
		Pick an interesting subset $\mathcal{S}(N)$ of the space $\M(N)$ of all
		$N \times N$ complex matrices.
		
		\item
		Put an interesting probability measure $\eta$ on $\mathcal{S}(N).$
		
		\item
		Select an interesting random variable (measurable function)
		$f:\mathcal{S}(N) \rightarrow \C.$
		
		\item
		Compute the expected value $\langle f \rangle$ of the random variable
		$f$ (possibly after rescaling) as a power series in $\frac{1}{N}:$
			\begin{equation}
				\langle f \rangle = \int_{\mathcal{S}(N)} f(M) \eta(dM) = 
				\sum_{g \geq 0} \frac{\varepsilon_g(f)}{N^g}.
			\end{equation}
			
		\item
		Realize that the coefficients $\varepsilon_g(f)$ occurring in the above perturbative
		expansion have an interesting combinatorial interpretation.
		
	\end{enumerate}
	\noindent
This informal discussion is meant to convey the impression that, from a combinatorial perspective,
matrix integrals may sometimes play the role of generating functions.  It sometimes happens
that a traditional generating function is difficult to obtain, but that by running the above
steps in reverse one can concoct a matrix integral which encodes a sequence of interest.  Furthermore, some matrix
models have special features that are very useful, and these features can be used to 
extract combinatorial information in the generating function spirit.  For example, it might be that an integral of 
interest can be exactly evaluated, thereby yielding an explicit generating function for the sequence
it encodes (see \cite{HZ} for a famous example).  Even if this is not the case, it often 
happens that matrix integrals interact well with more advanced analytical tools, such 
as orthogonal polynomials or integrable systems of differential equations (see e.g.
\cite{VM}), and this may again yield insight into combinatorial structure.

Let us present a matrix model for primitive factorizations.  For our space
of matrices we select the group $\U(N)$ of $N \times N$ complex unitary matrices.  Since $\U(N)$
is compact, it carries a unique left and right translation invariant probability measure,
the Haar measure $dU,$ which we take for our probability measure of interest.  Now we will 
select an interesting random variable, or rather class of random variables, $f:\U(N) \rightarrow \C.$
It would certainly be nice if we could compute the expected value $\langle P(U,U^{-1})\rangle$ of 
any polynomial function of the entries of $U$ and $U^{-1}=U^*,$ since we can approximate
a large class of functions on $\U(N)$ by polynomials in matrix coefficients (think 
Stone-Weierstrass/Peter-Weyl).  By linearity of the integral,
it suffices to consider the case where $P$ is a monomial.  Furthermore, an easy argument 
using the invariance of the Haar measure shows that the expected value of such a monomial
will be zero unless $P$ is of equal degree in the entries of $U$ and $U^*$ (think
integrals of the form $\int z^m \overline{z}^n dz$ over the unit circle $\U(1)$).  Thus 
we need only consider integrals of the form
	\begin{equation}
		\langle u_{i(1)j(1)} \overline{u}_{i'(1)j'(1)} \dots u_{i(n)j(n)} \overline{u}_{i'(n)j'(n)} \rangle
		= \int_{\U(N)} u_{i(1)j(1)} \overline{u}_{i'(1)j'(1)} \dots 
		 u_{i(n)j(n)}\overline{u}_{i'(n)j'(n)} dU,
	\end{equation}
where the lowercase $u_{ij}$'s are matrix elements and $i,j,i',j':\{1,\dots,n\} \rightarrow 
\{1,\dots,N\}$ are functions.  Integrals of this form are called \emph{$n$-point correlation 
functions} of matrix elements.  They are actually of considerable interest in 
mathematical physics \cite{DH,GT, Morozov} and free probability theory \cite{Collins}.  It is known that,
provided $N \geq n,$ the computation of the $n$-point functions can be reduced to 
the computation of ``permutation correlators''
	\begin{equation}
		\label{permCorrelators}
		\langle u_{11}\overline{u}_{1\pi(1)} \dots u_{nn}\overline{u}_{n\pi(n)} \rangle, \qquad
		\pi \in S(n).
	\end{equation}
Finally, recall that we have defined
	\begin{equation}
		a_{k,\mu} = \#\{\text{primitive factorizations of $\pi$ into $k$ transpositions}\},
	\end{equation}
and that this quantity can be non-zero only for $k$ of the form $k=n-\ell(\mu)+2g.$  Therefore
let us introduce the notation
	\begin{equation}
		\tilde{a}_{g,\mu}:=a_{n-\ell(\mu)+2g,\mu}.
	\end{equation}
It is not unreasonable to think of $\tilde{a}_{g,\mu}$ as a combinatorially motivated 
analogue of the usual Hurwitz number \cite{ELSV} $h_{g,\mu},$ obtained by replacing
the transitivity constraint with the primitivity constraint.

\begin{thm}
	\label{thm:matrixIntegral}
	Let $\mu$ be a partition of $n$ and let $\pi \in C_{\mu}$ be a permutation of 
	cycle type $\mu.$  Then, for any $N \geq n,$
		$$(-1)^{n-\ell(\mu)}N^{2n-\ell(\mu)}
		\langle u_{11}\overline{u}_{1\pi(1)} \dots u_{nn} \overline{u}_{n\pi(n)} \rangle
		= \sum_{g \geq 0} \frac{\a_{g,\mu}}{N^{2g}}.$$
\end{thm}

Unfortunately, we will not be able to say much about the proof of this result here.  Suffice 
to say that that Theorem \ref{thm:matrixIntegral} arises from two points of view regarding
the orthogonal projection of $\mathcal{M}(N)^{\otimes n}$ onto the commutant
	\begin{equation}
		\mathcal{C}_{\U(N)}(n) = \{T \in \M(N)^{\otimes n}: U^{\otimes n}T=TU^{\otimes n} \,\,\,
		\forall U \in \U(N)\}.
	\end{equation}
The first point of view involves the permutation correlators \eqref{permCorrelators} and 
the second involves understanding the element $(N+J_1)^{-1} \dots (N+J_n)^{-1}$ in 
the left-regular representation of $\C[S(n)];$ the equivalence of the two is, in a sense,
a manifestation of the Schur-Weyl duality between the representation theories of
$S(n)$ and $\U(N).$  We refer the interested reader to our articles \cite{MN,Novak} for 
further details regarding the proof of Theorem \ref{thm:matrixIntegral}, and its applications.

Let us conclude with the following comparison of unrestricted and primitive factorizations
of a full cycle: 
	\begin{align}
		h_{g,(n)} &= n^{n-2}n^{2g}{n-1+2g \choose n-1} \bigg{[} \frac{z^{2g}}{(2g)!} \bigg{]}
		\bigg{(} \frac{\sinh z/2}{z/2} \bigg{)}^{n-1} \\
		\tilde{a}_{g,(n)} &=\Cat_{n-1}{2n-2+2g \choose 2n-2}  \bigg{[} \frac{z^{2g}}{(2g)!} \bigg{]}
		\bigg{(} \frac{\sinh z/2}{z/2} \bigg{)}^{2n-2}.
	\end{align}
The first of these formulas is due to Jackson \cite{Jackson} (see also \cite{SSV}), while
the second is a consequence of Theorem \ref{thm:centralFactorial} together with Riordan's
exponential generating function for the central factorial numbers.
	
	
		
		
		


\end{document}